\documentclass{IEEEtran}
\usepackage{amsmath,amssymb,amsfonts}
\usepackage{algorithmic}
\usepackage{graphicx}
\usepackage{textcomp}
\def\BibTeX{{\rm B\kern-.05em{\sc i\kern-.025em b}\kern-.08em
    T\kern-.1667em\lower.7ex\hbox{E}\kern-.125emX}}

\usepackage[framed]{matlab-prettifier}

\usepackage{xcolor}

\usepackage{tikz}
\usepackage{circuitikz}
\usetikzlibrary{math}
\usepackage{mathrsfs}

\usepackage{siunitx}

\usepackage{placeins}

\ifCLASSOPTIONcompsoc
 \usepackage[caption=false,font=normalsize,labelfont=sf,textfont=sf]{subfig}
\else
 \usepackage[caption=false,font=footnotesize]{subfig}
\fi

\usepackage{pgfplots}
\pgfplotsset{compat=1.18} 
\usepgfplotslibrary{colorbrewer}
\usetikzlibrary{pgfplots.colorbrewer}
\usetikzlibrary{plotmarks}

\usepackage{cite}

\newcommand{\figref}[1]{Fig.~\ref{#1}}
\newcommand{\secref}[1]{Section~\ref{#1}}
\newcommand{\pM}{{(M)}}
\newcommand{\id}[0]{\mathbb{I}}
\newcommand{\eqcite}[2]{[\citenum{#1}, eq. (#2)]}
\newcommand{\rOne}[1]{#1}
\newcommand{\rTwo}[1]{#1}

\pgfplotsset{
  log x ticks with fixed point/.style={
      xticklabel={
        \pgfkeys{/pgf/fpu=true}
        \pgfmathparse{exp(\tick)}%
        \pgfmathprintnumber[fixed relative, precision=3]{\pgfmathresult}
        \pgfkeys{/pgf/fpu=false}
      }
  },
  log y ticks with fixed point/.style={
      yticklabel={
        \pgfkeys{/pgf/fpu=true}
        \pgfmathparse{exp(\tick)}%
        \pgfmathprintnumber[fixed relative, precision=3]{\pgfmathresult}
        \pgfkeys{/pgf/fpu=false}
      }
  }
}
    
\begin{document}
\title{Solver Performance of Accelerated \\ MoM for Connected Arrays}
\author{Harald Hultin, Lucas Åkerstedt, and B.L.G. Jonsson
\thanks{\rOne{This work was supported by the Swedish Foundation for
Strategic Research under Project ID20-0004.}}
\thanks{\rOne{Harald Hultin is with Saab Surveillance, 175 41 Järfälla, Sweden.}}
\thanks{\rOne{The authors are with KTH Royal Institute of Technology, EECS,
10044 Stockholm, Sweden (e-mail: haraldhu@kth.se).}}}

\maketitle

\begin{abstract}
Simulating and developing large rectangularly shaped arrays with equidistant interspacing is challenging as the computational complexity grows quickly with array size. However, the geometrical shape of the array, appropriately meshed, leads to a multilevel Toeplitz structure in the RWG-based Method of Moment impedance matrix representation that can be used to mitigate the increased complexity. This paper develops, presents and compares two different \rTwo{accelerated} solvers that both utilize \rTwo{the matrix} structure to determine antenna properties. Both methods use a novel mesh-partitioning algorithm and its associated \rTwo{data} representation, \rTwo{reducing storage and computational costs}. The first solver is an iterative method based on \rOne{multilevel fast Fourier transform} to accelerate matrix multiplications. The second solver approach is based on an extension of a fast direct Toeplitz solver, adapted to a block-matrix structure. \rTwo{This fast direct solver is demonstrated to have close to machine epsilon accuracy}. Both \rTwo{accelerated} methods are evaluated on two different array element types, \rOne{for arrays with up to 900 elements}. The results are compared with conventional direct and iterative matrix solvers. Improvements are seen in both the time and required storage to solve the problem. \rOne{The choice of the most efficient method depends on the residual thresholds in the iterative method, geometry of the element and frequency.} \rOne{Two different preconditioners for the iterative method are investigated to evaluate their performance. The two \rTwo{accelerated} methods vastly outperform regular matrix inversion methods.} 
\end{abstract}


\section{Introduction}
\label{sec:introduction}
\bstctlcite{IEEEexample:BSTcontrol}

\IEEEPARstart{M}{ethods} for simulating finite-size periodic structures are used in a range of applications across multiple disciplines. Such structures, composed of a unit cell repeated over a finite grid, include array antennas \cite{latha2021review}, reconfigurable intelligent surfaces \cite{liu2021reconfigurable}, and metasurfaces \cite{quevedo2019roadmap}.

This paper focuses on the simulation of array antennas, which are key components in a wide range of applications. These applications include defense \cite{spezio2002electronic, haupt2015antenna, talisa2016benefits}, communications \cite{ghosh2019inclusive, jiang2021road} and radio astronomy \cite{dewdney2009square, cavillot2019efficient} systems. A common trend across these areas is the transition to fully or partially digital arrays \cite{talisa2016benefits}. In a digital array, each element or sub-array incorporates radio-frequency (RF) amplification, followed by analog-to-digital (A/D) or digital-to-analog (D/A) conversion on receive or transmit, respectively. Time delays and amplitude weights are then applied in the digital domain, enabling flexible pattern control and multiple simultaneous beams \cite{mailloux2017phased}. These capabilities are critical for meeting growing demands on performance, adaptability, and multi-functionality \cite{talisa2016benefits, jiang2021road, harlakin2021compressive}. 

However, simulating large-scale or complex array antennas remains a computational challenge, particularly as system sizes and operational requirements continue to increase. Full-wave solvers, which are commonly used to model such systems, encounter scalability issues as the number of elements increases or the design complexity grows. One option is to use unit cell methods to approximate array performance, but to fully utilize the capabilities of the array antenna requires the correct behavior of each separate element in the array \cite{helander2017synthesis, emadeddin2024advancements, hultin2025wideband, jonsson2025farfield}

How to overcome the computational challenge of large problems is a well studied area, with several integral-based methods proposed\cite{ergul2014multilevel, He2022solving, kurz2002adaptive, zhao2005adaptive, maaskant2008fast, borm2003introduction, gujjula2022new, hackbusch2002blended, phillips1997precorrected, bleszynski1996aim, bleszynski2004block, kindt2003array, helander2017synthesis, brandt2024extended}. The multilevel fast multipole algorithm (MLFMA) \cite{ergul2014multilevel, He2022solving} evaluates weak interactions for a group of elements using a multipole representation. Another method to reduce the complexity is adaptive cross-approximation (ACA) \cite{kurz2002adaptive, zhao2005adaptive} that finds low-rank approximations for the weak interactions in an adaptive manner. Neither MLFMA nor ACA assumes that the impedance matrix contains any structure and therefore works well on general problems. ACA has also been shown to work well for antenna arrays in combination with \rTwo{the} characteristics basis function method \cite{maaskant2008fast}.

An alternative to using ACA is to use hierarchal matrices \cite{borm2003introduction}. In these matrices, the data is arranged in a binary tree structure and allows varying block sizes, with larger blocks away from the matrix diagonal represented by a low-rank approximation due to weaker inter-element interaction. Hierarchical matrix methods have been shown to accelerate the solution of integral equation problems \cite{borm2003introduction, gujjula2022new}, and have also been used with Toeplitz structures \cite{hackbusch2002blended}. 

\rOne{Another approach is to solve the integral equation problem on a uniformly spaced grid, such as the precorrected fast Fourier transform method (pFFT) \cite{phillips1997precorrected} and adaptive integral method (AIM) \cite{bleszynski1996aim}.} The periodicity of the grid then allows for a efficient solution of the problem via \rOne{fast Fourier transform (FFT)}. Similar FFT accelerated methods for unconnected arrays \cite{bleszynski2004block, kindt2003array, helander2017synthesis} utilize the periodicity of the array, which gives the impedance matrix a multilevel block-Toeplitz structure\rTwo{. This structure allows for a fast matrix-vector multiplication which enables the problem to be solved fast using an iterative solver, while also reducing memory usage}. 

However, the methods in \cite{bleszynski2004block, kindt2003array, helander2017synthesis} are only for unconnected arrays, the elements in the array must be separated by air. This excludes e.g. ground planes and most meta-surfaces. Brandt-M{\o}ller et al \cite{brandt2024extended} recently presented a method that works for connected arrays by using the discontinuous Galerkin method and high order, quadrilateral basis functions. In \cite{brandt2024extended}, it is mentioned that their method should work for RWG basis function \cite{rao1982electromagnetic}. 

In \cite{akerstedt2024partitioning} a novel decomposition method for array antennas with connected elements meshed with RWG basis functions is presented. That decomposition method is used in this paper to present a method to decompose and simulate antenna arrays efficiently. The method allows for the edges of array element meshes to overlap, assigning the current over the border to only one element without using a discontinuous Galerkin method. This comes with the feature that the border region of the array can be handled separately. These regions are not part of the Toeplitz structure. The full impedance matrix $Z$ has the form
\begin{equation}
\label{eq:Zblock}
    Z=\begin{bmatrix}
        Z_A & Z_B^T \\
        Z_B & Z_C
    \end{bmatrix},
\end{equation}
where $Z_A$ is the multilevel Toeplitz-structured matrix of the array elements, $Z_C$ is the impedance matrix of the border and $Z_B$ is the interaction between the elements and the border. 

In this paper, we describe and compare two different methods that use the structure \eqref{eq:Zblock} derived in \cite{akerstedt2024partitioning} to accelerate the array simulation. The first approach uses a multilevel FFT (MLFFT) to create a block-diagonal representation of $Z_A$, similar to \cite{bleszynski2004block, kindt2003array, helander2017synthesis, brandt2024extended}. This diagonalization allows for an accelerated multiplication in an iterative solver, \rOne{while storing only the unique block matrices} within $Z_A$. For better convergence in the iterative solver, we have implemented the preconditioner of \cite{kindt2003array} and a preconditioner similar to that of \cite{zhai2006analysis}. We also present a block variant of the Rybicki Toeplitz inversion given in \cite{press1986numerical}. This is a direct solver that can invert a 1-level block-Toeplitz matrix, in contrast to the MLFFT which can utilize all the levels of the multilevel structure of $Z_A$. Other similar methods may put requirements on the blocks in the block Toeplitz structure. For example \cite{watson1973algorithm} requires the block to be symmetric around both the primary and secondary diagonal. The Rybicki method presented in this paper works on \rTwo{invertible blocks with arbitrary structure}.

The performance of the \rTwo{two} presented methods is investigated by simulating two different array types of varying sizes and comparing the presented methods with conventional matrix inversions. The results show a definite performance improvement, both in time and \rTwo{scaling. Memory usage and especially memory scaling are improved by only saving the required blocks in $Z_A$}. For an array antenna with \rOne{$N_x$ columns and $N_y$ rows, the full size of $Z_A$ in \eqref{eq:Zblock} scales as $N_x^2N_y^2$. As described in \cite{akerstedt2024partitioning}, $Z_A$ for such a 2D array has an inherent multilevel block Toeplitz structure with two levels. On the top level the structure arises from identical interactions between rows of the array, and on the lower level from identical interactions between antenna elements within the rows.} By utilizing one level of the structure in the Rybicki method, the required data to \rTwo{store} $Z_A$ is reduced to scale as $N_x^2N_y$. Furthermore, the scaling is reduced even further by utilizing both levels in the structure in the MLFFT method, down to $N_xN_y$. Both MLFFT and Rybicki are faster and more efficient than the conventional methods. Whether MLFFT or Rybicki is the most efficient choice depends on the simulated structure and the desired relative residual of the MLFFT method.

Both the methods presented are explained explicitly, with code for the MLFFT given in App.~\ref{app:SLFFT}~and~\ref{app:DLFFT} and \rTwo{modifications} to the Rybicki method in \cite{press1986numerical} given in \secref{sec:Rybicki}. We present a thorough analysis of convergence and investigate how array size, different types of antenna elements, meshing and frequency effect the efficiency of the presented solvers. While MLFFT is faster for certain problems, Rybicki has a very stable performance over frequency for a given mesh and array size.

The paper is organized as follows. \secref{sec:theory} describes how a multilevel block-Toeplitz is block diagonalized, enabling a low complexity matrix-vector multiplication. \secref{sec:matinv} describes the methods used for matrix inversions. \secref{sec:simulation} describes how the methods in \secref{sec:matinv} are used to simulate an antenna array, and the methods used for comparison. In \secref{sec:results} the results of the performance comparison is presented and discussed. The paper ends with conclusions.

\section{Theory}
\label{sec:theory}
The MLFFT solver is accelerated by using a fast multiplication described in \secref{sec:toepmult}. This type of acceleration is used in several publications like \cite{brandt2024extended, bleszynski2004block, kindt2003array, helander2017synthesis}, but it is usually described in a rather condensed manner. In this section, we describe how the acceleration is implemented on a multilevel block-Toeplitz structure by starting with a scalar Toeplitz matrix, extend to the block-matrix case, and then introduce the one level and extend to the two level structure to give a thorough explanation of the theory behind the acceleration. This theory is used to derive the codes in \rTwo{Appendices}~\ref{app:SLFFT}~and~\ref{app:DLFFT}.
\subsection{Multilevel block-Toeplitz structure}
Multilevel block-Toeplitz structures, also called $L$-level block-Toeplitz structures \cite{choromanski2022block}, with $L$ levels have the following structure. On the highest level $L$, it is a block-Toeplitz structure:
\begin{equation}
\label{eq:Tblock}
    T^{[L]} = 
\begin{bmatrix}
R_{0} & R_{-1} & \cdots & R_{-N^{[L]}+1} \\
R_{1} & R_{0} &  & \vdots \\
\vdots &  & \ddots & R_{-1} \\
R_{N^{[L]}-1} & \cdots & R_{1} & R_{0} \\
\end{bmatrix}.
\end{equation}
If $L-1=l>0$, each block $R_n$ is a multilevel block-Toeplitz matrix with $l$ levels. If $l=0$, $R_n$ is here an unstructured matrix. \rTwo{The matrix $T^{[L]}$} is similar to the recursive block-Toeplitz matrix in \cite{lee1986fast}, though in a recursive block-Toeplitz matrix the lowest level consists of Toeplitz matrices rather than unstructured matrices. 

For finite periodic structures where each element has an identical mesh, the interaction between cells give this type of impedance matrix \cite{brandt2024extended, bleszynski2004block, kindt2003array, akerstedt2024partitioning}. The number of blocks on each level, $N^{[l]}$, is defined by the problem. For a one dimensional array with $N_x$ antenna elements, which in turn are discretized to $N_{e}$ mesh cells, the impedance matrix is a 1-level block-Toeplitz matrix with $N^{[1]}=N_x$ and $N^{[0]}=N_{e}$. Extending to a two dimensional array with $N_x\times N_y$ elements gives a 2-level block-Toeplitz matrix, with $N^{[2]}=N_y$, $N^{[1]}=N_x$ and $N^{[0]}=N_{e}$, which is the case that is used throughout this paper.

\subsection{Toeplitz Multiplication using DFT}
\label{sec:toepmult} 
A key property of a Toeplitz \rTwo{matrix}, and by extension a multilevel block-Toeplitz matrix, is the fast multiplication via FFT as shown in \cite{kindt2003array, brandt2024extended, alma99286799402456, choromanski2022block, lee1986fast}. This stems from the fact that a circulant matrix, $C$, is diagonalized by the Discrete Fourier Transform (DFT), $\mathcal{F}$:
\begin{equation}
    \mathcal{F} C \mathcal{F}^{-1} = \text{diag}(\mathcal{F}c_1),
\end{equation}
where $c_1$ is the vector consisting of the first column in $C$. Using $\mathcal{F}^{-1} \mathcal{F} = \id$, where $\id$ is the identity matrix, the product \rTwo{$Cu_c=v_c$} can be written as
\begin{equation}
    \label{eq:toeplitzmult}
    Cu_c =  \mathcal{F}^{-1} \mathcal{F} C \mathcal{F}^{-1} \mathcal{F}u_c = \mathcal{F}^{-1} \text{diag}(\mathcal{F}c_1) \mathcal{F}u_c.
\end{equation}

A Toeplitz matrix $T$ is not circulant, but can be padded to create the circulant matrix. Denote the circulant matrix created from $T$ by $C(T)$. To show the procedure, consider the following example. Take
\begin{equation}
\label{eq:T}
    T = 
\begin{bmatrix}
t_{0} & t_{-1} & t_{-2} \\
t_{1} & t_{0} & t_{-1} \\
t_{2} & t_{1} & t_{0} \\
\end{bmatrix}, \     T' = 
\begin{bmatrix}
t_{0} & t_{-1}  \\
t_{1} & t_{0}  \\
\end{bmatrix}.
\end{equation}
The circulant matrix $C(T)$ is populated:
\begin{equation}
\begin{split}
    C(T) = \begin{bmatrix}
        T & * \\
        * & T' \\
    \end{bmatrix},
\end{split}
\end{equation}
where $*$ are to be determined. Following the cyclic pattern in a circulant matrix, the missing entries follows as:
\begin{equation}
\label{eq:scalarC}
\begin{split}
    C(T) 
    =\begin{bmatrix}
    t_{0} & t_{-1} & t_{-2} & t_{2} & t_{1} \\
    t_{1} & t_{0} & t_{-1} & t_{-2} & t_{2} \\
    t_{2} & t_{1} & t_{0} & t_{-1} & t_{-2} \\
    t_{-2} & t_{2} & t_{1} & t_{0} & t_{-1}  \\
    t_{-1} & t_{-2} & t_{2} & t_{1} & t_{0}  \\
    \end{bmatrix}.
\end{split}
\end{equation}
The matrix multiplication $Tu=v$ can now be evaluated by
\begin{equation}
\label{eq:circmult}
    C(T) \begin{bmatrix}
        u \\
        0 \\
    \end{bmatrix}
    =
    \begin{bmatrix}
        Tu \\
        \gamma \\
    \end{bmatrix}
    =
    \begin{bmatrix}
        v \\
        \gamma \\
    \end{bmatrix}.
\end{equation}
Note that $0$ is here a zero vector of length $2$ and $\gamma$ a vector of unused values. To form $c_{1}$, the first column of $C$, take the first column of $T$. Then, take the elements in first row of $T$, except the first element, invert the order and append the the first column to create
\begin{equation}
\label{eq:cT1}
     c_{1}=
     \begin{bmatrix}
        t_0 & t_1 & t_2 & t_{-2} & t_{-1}
    \end{bmatrix}^T
    .
\end{equation}
That is, for a Toeplitz matrix of size $n\times n$, we only need to save $2n+1$ values to carry out the multiplication. Further, the multiplication is carried out with computational complexity in $\mathcal{O}(n\log(n))$ \rTwo{if the DFT} in \eqref{eq:toeplitzmult} \rTwo{is replaced by the FFT}.

\subsection{Multiplication of multilevel Toeplitz Matrix using DFT}
\label{sec:multitoepmult}
If the elements in $T$ and $C(T)$ in Section \ref{sec:toepmult} are replaced by block-matrices, it gives the single-level block matrices $T^{[1]}$ and $C^{[1]}(T^{[1]})$ which is defined analogously to \eqref{eq:scalarC}. The multiplication in \eqref{eq:circmult} is carried out in the same manner using the first block column of $C^{[1]}(T^{[1]})$ analogously to \eqref{eq:cT1}. \rTwo{The bolck column is a block vector, a block matrix with a single column of blocks}. However, $C^{[1]}(T^{[1]})$ is not diagonalized by the regular DFT. But as shown in, e.g. \cite{jelich2021efficient, kindt2003array}, it is block diagonalizable by a blockwise DFT. For a matrix structure of a 1D array with $N_x$ antenna elements, where each element is meshed with $N_e$ mesh cells, the 1-level block-wise DFT (1LDFT) $\mathcal{F}^{[1]}$ is defined as
\begin{equation}
\label{eq:F1}
\begin{split}
    \mathcal{F}^{[1]} &= \mathcal{F}_{N_x}\otimes \id_{N_e} 
    =      \begin{bmatrix}
        f_{1,1}\id_{N_e} & \hdots & f_{1,N_x}\id_{N_e} \\
        \vdots & \ddots & \vdots \\
        f_{N_x,1}\id_{N_e} & \hdots & f_{N_x,N_x}\id_{N_e} \\
    \end{bmatrix},
\end{split}
\end{equation}
where $\mathcal{F}_{N_x}$ is the DFT of size $N_x$, $\otimes$ is the Kroenecker product, $\id_{N_e}$ is the identity matrix of size $N_e$ and $f_{i,j}$ is an element of $\mathcal{F}_N$. The inverse 1LDFT is defined in the same manner, with $\mathcal{F}_N$ replaced by $\mathcal{F}_N^{-1}$.

With the 1LDFT, the multiplication \rTwo{$T^{[1]}u=v$} can be carried out in two steps. First, create the first block column $C^{[1]}_1(T^{[1]})$ analogously to Section \ref{sec:toepmult}:
\begin{equation}
\label{eq:circL1}
    C^{[1]}_1(T^{[1]}) =
     \begin{bmatrix}
        R_0^T &
        \hdots &
        R_{N_x-1}^T &
        R_{-N_x+1}^T &
        \hdots &
        R_{-1}^T 
    \end{bmatrix}^T,
\end{equation}
where $t_i$ in \eqref{eq:scalarC} is replaced by the corresponding submatrix $R_i$ in \eqref{eq:Tblock}. As each submatrix $R_i$ has dimension $N_e\times N_e$, the matrix $C^{[1]}_1(T^{[1]})$ has dimension $(2N_x-1)N_e \times N_e$.

Second, a block analogue to \eqref{eq:toeplitzmult} is required. For a general level $L$, the product can be written as
\begin{multline}   
\label{eq:multL}
    C^{[L]}(T^{[L]})u^{[L]} = \\ \left(\mathcal{F}^{[L]}\right)^{-1} \text{diag}_{N^{[0]}}\left(\mathcal{F}^{[L]}C^{[L]}_1(T^{[L]})\right)\mathcal{F}^{[L]}u^{[L]}=v^{[L]} \\
    u^{[L]},v^{[L]} \in \mathbb{C}^{N_C},
\end{multline}
where $N_C=N^{[0]}\prod_{i=1}^L(2N^{[i]}-1)$ is the dimension of $C^{[L]}$, $\text{diag}_{N^{[0]}}$ denotes a block-diagonal matrix with block size $N^{[0]}$ and vectors $u^{[L]}$ and $v^{[L]}$ are padded versions of $u$ and $v$. For the \rTwo{one dimensional array, $L=1$}, the vectors are defined as
\begin{equation}
\label{eq:pad1}
    u^{[1]} = \begin{bmatrix}
        u \\
        0 \\
    \end{bmatrix}, \:
        v^{[1]} = \begin{bmatrix}
        v \\
        \gamma \\
    \end{bmatrix},
\end{equation}
where the desired product $v$ lies in the first $N^{[1]}N^{[0]}$ elements of \rTwo{$v^{[1]}$. Using \eqref{eq:multL} and \eqref{eq:pad1} the matrix multiplication for $L=1$} can be carried out by only storing $(2N^{[1]}-1)\left(N^{[0]}\right)^2$ elements of the matrix.

As shown, operating on a 1-level block-Toeplitz structure is largely analogous to the scalar case. For a 2-level block-Toeplitz structure, some more consideration is required. The 2-level block-wise DFT (2LDFT) $\mathcal{F}^{[2]}$ is defined similarly to \rTwo{the 1LDFT} \eqref{eq:F1}:
\begin{equation}
\label{eq:F2}
    \mathcal{F}^{[2]} = \mathcal{F}_{N^{[2]}} \otimes \left(\mathcal{F}_{N^{[1]}} \otimes I_{N^{[0]}}\right) = \mathcal{F}_{N^{[2]}} \otimes \mathcal{F}^{[1]}
\end{equation}
where $N^{[2]}$ is the number of blocks on level 2. To create $C^{[2]}(T^{[2]})$ from $T^{[2]}$, the blocks need to be rearranged. The Toeplitz matrix of level 2, $T^{[2]}$, is structured as follows:
\begin{equation}
    T^{[2]} = 
    \begin{bmatrix}
        T_0^{[1]} & \hdots & T_{-N^{[1]}+1}^{[1]} \\
        \vdots & \ddots & \vdots \\
        T_{N^{[1]}-1}^{[1]} & \hdots & T_0^{[1]} \\
    \end{bmatrix},
\end{equation}
where $T_n^{[1]}$ are 1-level block-Toeplitz matrices.

This structure can, similar to creating $C^{[1]}(T^{[1]})$, be used to create the block diagonalizable matrix $C^{[2]}(T^{[2]})$:
\begin{multline}
    C^{[2]}(T^{[2]}) = \\
    \begin{bmatrix}
        C^{[1]}(T_0^{[1]}) & C^{[1]}(T_{-1}^{[1]}) & \hdots & C^{[1]}(T_1^{[1]}) \\
        \vdots & & & \vdots \\
        C^{[1]}(T_{N^{[1]}-1}^{[1]}) & \ddots & & \vdots \\
        C^{[1]}(T_{-N^{[1]}+1}^{[1]}) & & \ddots & \vdots \\
        \vdots & & & \vdots \\
        C^{[1]}(T_{-1}^{[1]}) & \hdots & \hdots&C^{[1]}(T_0^{[1]}) \\
    \end{bmatrix},
\end{multline}
which has a structure similar to \eqref{eq:scalarC} with $t_i$ replaced by $C^{[1]}(T_i^{[1]})$. The first block column of $C^{[2]}(T^{[2]})$, with blocks of size $N^{[0]}\times N^{[0]}$, is found as follows. For each $T_n^{[1]},\: n=1,\dots, N^{[1]}, -N^{[1]}, \dots, -1$, form a block column $C_1^{[1]}(T_n^{[1]})$ as in \eqref{eq:circL1}. Construct the complete block column $C_1^{[2]}(T^{[2]})$ as
\begin{multline}
\label{eq:circL2}
    C_1^{[2]}(T^{[2]}) = 
        \left[\begin{matrix}
            C_1^{[1]}(T_0^{[1]})^T &
            \hdots &
            C_1^{[1]}(T_{N^{[1]}-1}^{[1]})^T
        \end{matrix}\right. \\
        \left. \begin{matrix}
            C_1^{[1]}(T_{-N^{[1]}+1}^{[1]})^T &
            \hdots &
            C_1^{[1]}(T_{-1}^{[1]})^T
        \end{matrix} \right]^T,
\end{multline}
which gives the same structure as in \cite{brandt2024extended}, Fig. 1. From this structure it can also be seen that the zero-padding of $u$ and $v$ becomes more complicated, also shown in \cite{brandt2024extended}. Take the original vector $u$ and divide it into $N^{[2]}$ subvectors of length $N^{[1]}N^{[0]}$:
\begin{equation}
    u = \begin{bmatrix}
        u_0^T &
        u_1^T &
        \hdots &
        u_{N^{[2]}-1}^T
    \end{bmatrix}^T.
\end{equation}
To perform the multiplication $C^{[2]}(T^{[2]})u^{[2]}=v^{[2]}$ correctly, the vector needs to be padded as in \eqref{eq:pad1}, but first block wise on the first level and then padded on the second level:
\begin{multline}
\label{eq:xc2}
    u^{[2]} = 
    \left[\begin{matrix}
        u_0^T &
        0_{(N^{[1]}-1)N^{[0]}} &
        u_1^T &
        0_{(N^{[1]}-1)N^{[0]}} &
        \hdots
    \end{matrix}\right. \\
    \left. \begin{matrix}
        u_{N^{[2]}-1}^T &
        0_{(N^{[1]}-1)N^{[0]}} &
        0_{(N^{[2]}-1)(2N^{[1]}-1)N^{[0]}}
    \end{matrix} \right]^T,
\end{multline}
where $0_N$ is a zero row-vector of length $N$. The multiplication can now be performed using \eqref{eq:multL} with $L=2$. The resulting $v^{[2]}$ has a similar structure to \eqref{eq:xc2}, but with unwanted data instead of zero vectors. Finally, $v_n$ can be extracted from positions $(1, \dots, N^{[1]}N^{[0]}) + n(2N^{[1]}-1)N^{[0]}$ of $v^{[2]}$ giving
\begin{equation}
    v = \begin{bmatrix}
        v_0^T &
        v_1^T &
        \hdots &
        v_{N^{[2]}-1}^T 
    \end{bmatrix}^T.
\end{equation}
\section{Method}
\label{sec:matinv}
Method of Moments (MoM) gives a typical linear equation system on the form $ZI=V$, where $Z$ is an impedance matrix, $I$ is the unknown current and $V$ the excitation voltage. \rTwo{In this paper}, an in-house MoM code is used. The code method to partition the matrix into a \rTwo{multilevel block} Toeplitz structure is described in \cite{akerstedt2024partitioning}. To find the unknown current, two different methods are presented in this work: \rTwo{one using the multiplication from Section \ref{sec:theory} accelerated by multilevel FFT (MLFFT)} in an iterative solver, and \rTwo{one using} a block generalization of the \rTwo{Rybicki method \cite{press1986numerical}}.

\subsection{Multilevel FFT}
\label{sec:MLFFT}
To efficiently implement the ideas in \ref{sec:multitoepmult} and enable an efficient matrix multiplication, 1LDFT and 2LDFT are replaced with MLFFT. Here, it is implemented in Matlab by performing an FFT over indices in the blocks. For a 1-level block-Toeplitz structure, it can be implemented as in Appendix \ref{app:SLFFT}. The key observation here is that taking each $N^{[0]}$:th row of \eqref{eq:F1} and multiplying it with each $N^{[0]}$:th row of a vector $u$ just gives the DFT of those elements in $u$, which can be saved in the corresponding rows of the output. This can be done for all sets of rows, giving the block-wise FFT. The block-wise FFT is applied column-wise and consists of $N^{[0]}$ FFTs of length $(2N^{[1]}-1)=L_F$ and is applied on a block matrix on the form of \eqref{eq:circL1} with length $L_FN^{[0]}$ and width $W$. Each FFT here has a computational complexity of $\mathcal{O}(L_F\log (L_F))$. An FFT is applied for each of the $N^{[0]}$ set of rows, once per column over a total of $W$ columns, resulting in a total computational complexity of $\mathcal{O}(WN^{[0]}N^{[1]}\log (N^{[1]}))$.

In a similar manner to defining a MLFFT for the 1-level block-Toeplitz structure, a MLFFT for a 2-level structure can be implemented considering two block indices instead of one as shown in Appendix \ref{app:DLFFT}. For each of the $W$ columns, first apply $N^{[2]}N^{[0]}$ FFTs of length $N^{[1]}$, with complexity $\mathcal{O}(WN^{[2]}N^{[0]}N^{[1]}\log (N^{[1]}))$. Second, apply $N^{[1]}N^{[0]}$ FFTs of length $N^{[2]}$, complexity $\mathcal{O}(WN^{[1]}N^{[0]}N^{[2]}\log (N^{[2]}))$. The first operation dominates if $N^{[1]}>N^{[2]}$ and vice versa.

\subsection{Preconditioner for \rTwo{MLFFT}}
\label{sec:precond}
Consider the array \rTwo{antenna} described in Section~\ref{sec:introduction}. Due to the \rOne{weakly} singular nature of the Green's function, interactions over short electrical distances become very strong and interactions over long electrical distances become weak, leading to $Z$ being an ill-conditioned matrix. \rTwo{Strong interactions are largely concentrated to the self interaction of elements and to neighboring elements, though this depends element geometry.} The condition number of the problem can be improved by including a preconditioner in the matrix equation. In this paper, we use a left preconditioning strategy with the preconditioner $P$, rewriting the problem as
\begin{equation}
    P^{-1}ZI = P^{-1}V,
\end{equation}
where the goal is for $P^{-1}Z$ to be better conditioned compared to $Z$. 

One approach is to use an incomplete LU-decomposition (ILU) \cite{saad2003iterative} on $Z$ which gives a decomposition such that $ Z \approx LU = P$. Then, the preconditioner $P^{-1} = \left(LU\right)^{-1}$ can be applied efficiently using back substitution. However, even if an incomplete LU of $Z$ is cheaper than a regular LU, both $L$ and $U$ have the same size as $Z$. A way to address this is to use the MATLAB function $ilu$ as it outputs the two matrices in a sparse matrix format to reduce the memory footprint. Even so, this approach yields a preconditioner that, to a large extent, counteracts the memory savings of \eqref{eq:multL}.

A more memory efficient preconditioner has been demonstrated by Kindt et al. \cite{kindt2003array}. This preconditioner $P'_K$ is a block-diagonal matrix, where each block is the self-interaction of one element in the matrix $Z_{n,n}$, which has already been computed as it is the block-diagonal part of $Z_A$ in \eqref{eq:Zblock} on the lowest level. As the matrix is block diagonal, it can be inverted block wise. Further, as each block is identical, one LU-decomposition of $Z_{n,n}$ is sufficient to apply the inverse of $P'_K$ block wise. To cite \cite{saad2003iterative}: "For example, preconditioners can be derived from knowledge of the original physical problems from which the linear system arises". With the preconditioner $P'_K$, we use the knowledge that $Z_A$ is largely dominated by the strong block diagonal self-interaction. The matrix $P'_K$ has $N_{e}^2$ elements.

This knowledge of the underlying problem can be extended to another block-diagonal preconditioner $P'_Z$, similar to that of Zhai et al. \cite{zhai2006analysis}. A natural choice of sub-array in our method is one row of the matrix. Such a preconditioner consists of the block diagonal of $Z_A$, but on one level higher in the multilevel block-Toeplitz structure. Each block $Z^{[1]}_{n,n}$ is a block-Toeplitz matrix which describes the self interaction of a row of elements in the array. The preconditioner $P'_Z$ can be applied in the same manner as $P'_K$, though the block is larger with $N_{x}^2N_{e}^2$ elements.

The preconditioners $P'_K$ and $P'_Z$ are defined for $Z_A$ in \eqref{eq:Zblock}. However, $Z_A$ is only one, albeit the largest, part of $Z$. To account for the margin parts of $Z$, the full preconditioners $P_K$ and $P_Z$ are defined as
\begin{equation}
    P_X=\begin{bmatrix}
        P'_X & 0 \\
        0 & Z_C
    \end{bmatrix},
\end{equation}
where the index $X$ is either $K$ or $Z$. As the self interaction matrix of the border $Z_C$ is small, its full LU-decomposition is used as a preconditioner for the border. The LU-decomposition of $Z_C$ is applied in the same block-wise manner as above.

\subsection{Iterative Inverse with Multiple Right Hand Sides}
\label{sec:vectorization}
Given the fast multiplication procedure from Section \ref{sec:theory}, the MLFFT in Section \ref{sec:MLFFT} and the preconditioner in Section \ref{sec:precond}, the linear equation system $ZI_i=V_i$ can be solved using an iterative method. Here $V_i$ is an excitation vector corresponding to applying the desired excitation voltage to antenna port $i$. This excitation vector is passed to GMRES in MATLAB together with a function for the multiplication and a function for applying the preconditioner in order to solve for $I_i$, the excited current. This current can in turn be used to find the excited fields.

In addition to being a computationally complex problem, array antennas are often solved with the goal of the embedded element patterns (EEPs), i.e. to solve the matrix equation for every port excitation $V_i$. A straightforward way to achieve this is to run the solver once for each element, but that scales with computational complexity of $\mathcal{O}(N^{[2]}N^{[1]})$. 

Instead of solving each element one by one, consider the problem on matrix form:
\begin{equation}
\begin{split}
    \label{eq:ZIVmat}
    Z
    \begin{bmatrix}
        I_1 & \hdots & I_{N^{[2]}N^{[1]}}
    \end{bmatrix}
    &=
    \begin{bmatrix}
        V_1 & \hdots & V_{N^{[2]}N^{[1]}}
    \end{bmatrix}  \\
    = ZI_M &= V_M.
\end{split}
\end{equation}
The MLFFT based multiplication works for multiplying with a matrix, but standard GMRES only handles vector multiplicands and vector right hand sides. As shown in \cite{simoncini1995iterative}, there are methods to extend the underlying Krylov subspace to represent matrices rather than vectors. But rather than using a rewritten GMRES, the global Krylov subspace may be achieved using vectorization. As shown in \cite{alma99286799402456}, \eqref{eq:ZIVmat} is equivalent to
\begin{equation}
\label{eq:ZIVvec}
    (\id \otimes Z) \text{vec}(I_M) = \text{vec}(V_M),
\end{equation}
where $\text{vec}()$ is the vectorization operator \cite{alma99286799402456}. In MATLAB, the vectorization operator corresponds to reshaping to a matrix with one column, and the inverse by reshaping back to the original dimensions. This takes virtually no computational effort, as this reshaping only changes how MATLAB reads the variable from memory. As \eqref{eq:ZIVmat} and \eqref{eq:ZIVvec} are equivalent, the matrix multiplications can be performed as in \eqref{eq:ZIVmat}. There is no need to explicitly form the huge matrix $\id \otimes Z$. Then the vectorized right hand side (RHS) of \eqref{eq:ZIVvec} can be passed to GMRES. The Krylov subspace created by GMRES will then consists of the vectorized matrices in the subspace of the global GMRES described in \cite{simoncini1995iterative}.

\subsection{The Rybicki Method}
\label{sec:Rybicki}
The Rybicki method for a regular non-symmetric Toeplitz matrix is detailed in \cite{press1986numerical}. The Rybicki method is a bordering method \cite{ems2024bordering}, a recursive procedure that solves $T^{[1]}x=y$ through a sequence of solutions $x^{(M)}$ that converge in the last step when $M=N$. In our case $R_n$ is a matrix rather than a scalar as in \cite{press1986numerical}. \rOne{We update the algorithm in \cite{press1986numerical} with respect to matrices. The underlying equations \eqref{xmpo}~to~\eqref{eq:GH} corresponds to the key equations in \cite[Ch.~2.8]{press1986numerical}}.

First, \eqcite{press1986numerical}{2.8.19} it is rewritten as
\begin{multline}
    \label{xmpo}
	x_{M+1}^{(M+1)} = \left(\sum_{j=1}^M R_{M+1-j}G_{M+1-j}^\pM -R_0\right)^{-1} \\ \left(
	\sum_{j=1}^M R_{M+1-j}x_j^\pM-y_{M+1}\right). 
\end{multline}
Similarly, rewrite \eqcite{press1986numerical}{2.8.18} as
\begin{equation}\label{xi}
 x_j^{(M+1)} = x_j^\pM - G_{M+1-j}^\pM x_{M+1}^{(M+1)},
\end{equation}
\eqcite{press1986numerical}{2.8.24} as
\begin{multline}
    \label{eq:GM1}
	G_{M+1}^{(M+1)} = \left(  \sum_{j=1}^M R_{j-(M+1)}H_{M+1-j}^\pM-R_0
	\right)^{-1} \\ \left(\sum_{j=1}^M R_{j-(M+1)}G_j^\pM - R_{-(M+1)}\right)
\end{multline}
and \eqcite{press1986numerical}{2.8.23} as
\begin{multline}
\label{eq:HM1}
	H_{M+1}^{(M+1)} = 
	\left(\sum_{j=1}^M R_{M+1-j}G_{M+1-j}^\pM -R_0\right)^{-1} \\ \left(
	\sum_{j=1}^M R_{M+1-j}H_j^\pM-R_{M+1}\right). 
\end{multline}

When considering \eqcite{press1986numerical}{2.8.27}, the order in the multiplication is changed to take the non-comutability of matrices in to account by
\begin{equation}\label{HM}
H_{M+1-j}^{(M+1)} = H_{M+1-j}^{(M)} - G^{(M)}_{j}H_{M+1}^{(M+1)},\ j=1,\ldots M.
\end{equation}
Similarly, the order in \eqcite{press1986numerical}{2.8.25} is changed to
\begin{equation}
\label{eq:GH}
\begin{aligned}
    G_{j}^{(M+1)}= G_{j}^\pM- H_{M+1-j}^\pM G_{M+1}^{M+1}, \\
    H_{j}^{(M+1)}= H_{j}^\pM- G_{M+1-j}^\pM H_{M+1}^{M+1}.
\end{aligned}
\end{equation}

We observe that the obtained equations are almost identical to the scalar version, except for that iterated variable and the generator has exchanged place to respect the non-comutability of the matrices. \rOne{With these modifications, the subroutine given by \cite{press1986numerical} can be \rTwo{rewritten} to handle matrices with small changes to the code.}

This rewritten Rybicki method works for a block-Toeplitz structure, but has not been modified to fully utilize the multilevel structure of $Z_A$. Here, this means that the Rybicki method is implemented on level 2, where there is a block-Toeplitz structure and each block, i.e. each element in \textbf{R}, is a full block-Toeplitz matrix. Even so, this reduces the memory usage from the full matrix $Z_A$ with $N_y^2N_x^2N_{e}^2$ elements down to a structure with $(2N_y-1)N_x^2N_{e}^2$ elements. Note that this method works for any set of invertible $R_n$, unlike e.g. \cite{watson1973algorithm} \rOne{that has additional requirements} on symmetry.

\subsection{Including non-Toeplitz Blocks}
\label{sec:method_non_toep}
As mentioned in Section \ref{sec:introduction}, the impedance matrix given by \cite{akerstedt2024partitioning} contains non-Toeplitz structured blocks from the border of the array. Therefore, the two methods described in this section must be modified to include these blocks.
\subsubsection{MLFFT}
The MLFFT method is based on matrix multiplication, which can be performed blockwise. For a single RHS, partition the current vector $I=[I_A^T \ I_C^T]^T$ and the excitation vector $V=[V_A^T \ V_C^T]^T$ into the values on the array elements (denoted by $A$) and the values on the border (denoted by $C$). Then, multiply blockwise:

\begin{equation}
\label{eq:ZIVblock}
    ZI =
    \begin{bmatrix}
        Z_A & Z_B^T \\
        Z_B & Z_C
    \end{bmatrix}
    \begin{bmatrix}
        I_A \\
        I_C
    \end{bmatrix} 
    =
    \begin{bmatrix}
        Z_AI_A + Z_B^TI_C \\
        Z_BI_A + Z_CI_C \\
    \end{bmatrix}
    =
    \begin{bmatrix}
        V_A \\
        V_C
    \end{bmatrix}.
\end{equation}
Multiplying with $Z_A$ is accelerated as in Section \ref{sec:theory}, while the smaller margin matrices $Z_B$ and $Z_C$ are multiplied using matrix multiplication.

\subsubsection{Rybicki}
The Rybicki method only works on a Toeplitz structure, i.e. the problem must be on the form $T^{[1]}x=y$. Here, Schur decomposition is used to rewrite the matrix equation to only invert $Z_A$, moving the border matrices $Z_B$ and $Z_C$ the right hand side. Once this matrix equation has been solved efficiently, the solution of the original problem can be found by a relatively small inversion.

Given the first row of blocks in \eqref{eq:ZIVblock}, $I_A$ can be written as
\begin{equation}
    I_A = Z_A^{-1}(V_A - Z_B^TI_C).
\end{equation}
Inserting this in the second row of blocks gives
\begin{equation}
    Z_B Z_A^{-1} V_A - Z_B Z_A^{-1} Z_B^TI_C + Z_C I_C = V_C.
\end{equation}
Isolating $I_C$ gives
\begin{equation}
    (Z_C-Z_B Z_A^{-1} Z_B^T) I_C = V_C-Z_B Z_A^{-1} V_A.
\end{equation}
Let $U=Z_A^{-1} V_A$ and $F=Z_A^{-1} Z_B^T$ and solve the problem
\begin{equation}
    Z_A \begin{bmatrix}
        U & F
    \end{bmatrix}
    = \begin{bmatrix}
        V_A & B^T
    \end{bmatrix}
\end{equation}
for 
$\begin{bmatrix}
    U & F
\end{bmatrix}$.
The currents can then be found using
\begin{equation}
    I_C = (Z_C-Z_BF)^{-1}(V_C-Z_BU),
\end{equation}
a relatively small inversion, and
\begin{equation}
    I_A = U-F I_C.
\end{equation}

\section{Array Simulation}
\label{sec:simulation}
\begin{figure}[t]
\centering
\subfloat[]{
\resizebox{0.35\columnwidth}{!}{%
\begin{tikzpicture}[scale=0.25, font=\huge]
\draw[very thick] (\lf1/2,\h1+\t1/2) -- (\lf1/2+\l1,\h1+\t1/2) -- (\lf1/2+\l1+\lt1,\h1+\t2/2) -- (\w1/2,\h1+\t2/2) -- (\w1/2,\h1-\t2/2) -- (\lf1/2+\l1+\lt1,\h1-\t2/2) -- (\lf1/2+\l1,\h1-\t1/2) -- (\lf1/2,\h1-\t1/2) -- (\lf1/2,\h1+\t1/2);

\draw[<->] (\w1/2+1,0) -- node[right] {$h$} (\w1/2+1,\hd1-0.3);

\draw[<->] (-15, 1) -- node[above] {$w$} (15, 1);

\draw[very thick] (-15, 0) --(15, 0);

\draw[blue, dashed, very thick] (\w1/2,\h1+\t2/2) -- (\w1/2,\h1-\t2/2) -- (\w1/2-\c1,\h1-\t2/2) -- (\w1/2-\c1,\h1+\t2/2) -- (\w1/2,\h1+\t2/2);

\foreach \pos in {-7,...,6} {
    \node at (2*\pos,-2) { \tikz\draw [very thick] (-0,-0) -- (-1,-1);} ;}

\draw[red, fill] (-\lf1/2,\h1+\t1/2) -- (\lf1/2,\h1+\t1/2) -- (\lf1/2,\h1-\t1/2) -- (-\lf1/2,\h1-\t1/2) -- (-\lf1/2,\h1+\t1/2) ;

\begin{scope}[yscale=1,xscale=-1]
    \draw[very thick] (\w1/2,\h1+\t2/2) -- (\w1/2,\h1-\t2/2) -- (\w1/2-\c1,\h1-\t2/2) -- (\w1/2-\c1,\h1+\t2/2) -- (\w1/2,\h1+\t2/2);
  \draw[blue, dashed, very thick] (\lf1/2,\h1+\t1/2) -- (\lf1/2+\l1,\h1+\t1/2) -- (\lf1/2+\l1+\lt1,\h1+\t2/2) -- (\w1/2,\h1+\t2/2) -- (\w1/2,\h1-\t2/2) -- (\lf1/2+\l1+\lt1,\h1-\t2/2) -- (\lf1/2+\l1,\h1-\t1/2) -- (\lf1/2,\h1-\t1/2) -- (\lf1/2,\h1+\t1/2); % Mirror Image
  
\end{scope}
\end{tikzpicture}
}\label{fig:geometry_side}}
\subfloat[]{
\resizebox{0.45\columnwidth}{!}{%
    \begin{tikzpicture}[scale=4, thick]
    
    \draw[dotted] (-0.55,-0.55) -- (-0.55,0.55) -- (0.55,0.55) -- (0.55,-0.55) -- (-0.55,-0.55);

    \draw[thick] (-0.5,-0.375) -- (-0.5,0.375) -- (0,0.01) -- (0.5,0.375) -- (0.5,-0.375) -- (0,-0.01) -- (-0.5,-0.375);
    
    \draw[Set1-A, ultra thick] (0,-0.01) -- (0,0.01);

    \draw[<->] (-0.5,0.4) -- node[below] {$0.5\lambda$}(0.5,0.4);
    \draw[<->] (-0.55,-0.5) -- node[above] {$0.55\lambda$}(0.55,-0.5);
    \draw[<->] (0.55,0.375) -- node[above, rotate=270] {$W$}(0.55,-0.375);
    \draw[<->] (0.7,-0.55) -- node[above, rotate=270] {$0.55\lambda$}(0.7,0.55);

    \end{tikzpicture}
    }\label{fig:bowtie_geo}
    
}
\caption{Antenna unit cell geometries, feed shown in red. (a) CSA unit cell from the side, $h=\qty{20.7}{\milli\meter}$ and $w=\qty{30}{\milli\meter}$. The unit cell is \qtyproduct{30 x 30}{\milli\meter} seen from above. (b) bowtie unit cell from above, the bowtie is placed $0.25\lambda$ over a ground plane.} 
\label{fig:CSA}
\end{figure}
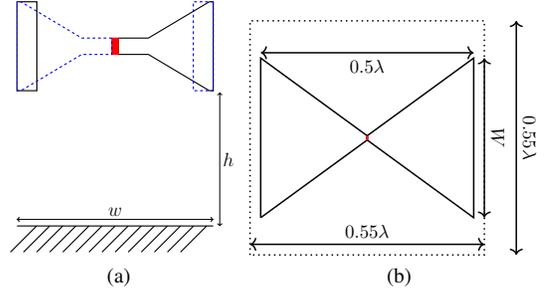

To evaluate the performance of the different matrix inversion methods in Section~\ref{sec:matinv}, two different types of array antennas are simulated: the Current Sheet Array (CSA) element of \cite{hultin2025wideband} illustrated in \figref{fig:geometry_side} and the bowtie element illustrated in \figref{fig:bowtie_geo}. The CSA element is described in full in \cite{hultin2025wideband} including a full geometrical description. Here, the CSA is simulated at a low frequency of \qty{1}{\giga\hertz} and a high frequency of \qty{4}{\giga\hertz}, with $462$ unknowns per element. The bowtie element has two variants. A wider, densely meshed element with $W=0.375\lambda$ and $1081$ unknowns per element is simulated at \qty{150}{\mega\hertz} and \qty{225}{\mega\hertz}, and a narrower element with $W=0.25\lambda$ and $253$ unknowns per element simulated at \qty{150}{\mega\hertz}.

For $Z$ such that the full matrix can be saved in memory, the methods in Section~\ref{sec:matinv} are compared to inversion using \emph{mldivide}, and using GMRES with a $P_K$ preconditioner. Note that the MLFFT methods in Section~\ref{sec:matinv} also use GMRES, but with the here proposed more efficient MLFFT matrix multiplication.

\tikzmath{
\t1 =2.5; \t2 = 13.6; \l1 = 4.59; \lt1 = 9.31; \c1 = 3.04; \h1 = 27.5; \w1 = 30; \lf1 = 1.0; \d1 = 1.0; \hd1 = \h1-\t2/2;
}

Here, a simulation consists of:
\begin{enumerate}
    \item read the mesh file of the antenna element,
    \item define simulation parameters, e.g. array size and frequency,
    \item create the unique sub-matrices of $Z$ in \eqref{eq:Zblock} as described in \cite{akerstedt2024partitioning},
    \item (conditional) fill in a full matrix variant of $Z$ (for \emph{mldivide} or GMRES) or fill in the lowest Toeplitz level of $Z$ (for the Rybicki method),
    \item solve $ZI=V$.
\end{enumerate}
\rOne{All simulations are run on a computer with an Intel Xeon Gold 5217 processor and 191 GB RAM.}

By default in this publication, $V$ is a matrix where each column $c$ consists of zeros, except for the feeding edge of array element $c$ which is set to $1$. The corresponding column $c$ of $I$ then contains the currents of the array from exciting the feeding point of element $c$. By solving one element at a time, arbitrary array excitations can be calculated in post-processing via super-position.

\section{Results}

\label{sec:results}

In order to investigate the performance of the presented methods, various array sizes with the different unit cell geometries were tested. The main metrics of interest is how fast the solvers are and how they scale with array size. The time to create the sub-matrices of $Z$ is given as "Con." for construction. The reference solutions solved using \emph{mldivide} and GMRES are shortened as "mld" and "G", their solution times include filling in the full matrix with the data from the construction step. 

The Rybicki solver, here shortened as "R", also includes the time to fill in the matrices on level 1, whereas the MLFFT solvers work with no fill-in. The MLFFT solver is investigated using three different configurations: M1, M2 and M3. in \secref{sec:precond}, two different preconditioners are defined; the smaller $P_K$ and the larger $P_Z$. Configuration M1 uses a $P_K$ preconditioner with a vectorized RHS, M2 a $P_Z$ preconditioner with vectorized RHS and M3 a $P_K$ preconditioner that solves the RHS column by column. The vectorization of the RHS is described in \secref{sec:vectorization}.

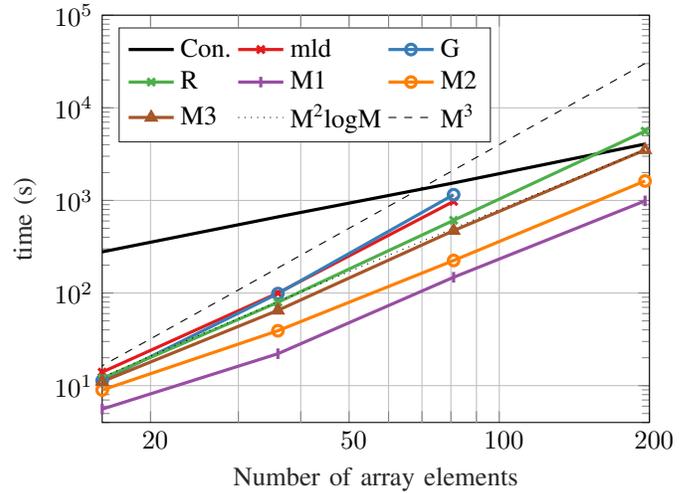
\begin{figure}[t]
    \centering
    % This file was created by matlab2tikz.
%
%The latest updates can be retrieved from
%  http://www.mathworks.com/matlabcentral/fileexchange/22022-matlab2tikz-matlab2tikz
%where you can also make suggestions and rate matlab2tikz.
%
\colorlet{mycolor1}{black}%
\colorlet{mycolor2}{Set1-A}%
\colorlet{mycolor3}{Set1-B}%
\colorlet{mycolor4}{Set1-C}%
\colorlet{mycolor5}{Set1-D}%
\colorlet{mycolor6}{Set1-E}%
\colorlet{mycoloryellow}{Set1-F}%
\colorlet{mycolor7}{Set1-G}%
\begin{tikzpicture}

\begin{axis}[%
width=\columnwidth,
height=7cm,
xmode=log,
xmin=16,
xmax=200,
xminorticks=true,
xlabel style={font=\color{white!15!black}},
xlabel={Number of array elements},
ymode=log,
ymin=4,
ymax=100000,
yminorticks=true,
xmajorgrids,
xminorgrids,
ymajorgrids,
ylabel style={font=\color{white!15!black}},
ylabel={time (s)},
xtick={10,20, 50, 100, 200},
log x ticks with fixed point,
extra x ticks={30, 40, 60, 70, 80, 90},
extra x tick label=\empty,
axis background/.style={fill=white},
legend style={at={(0.03,0.97)}, anchor=north west, legend cell align=left, align=left, draw=white!15!black, legend columns=3}
]
\addplot [color=mycolor1, very thick]
  table[row sep=crcr]{%
16	278.232731\\
36	664.3995005\\
81	1544.2765046\\
196	4060.4189562\\
};
\addlegendentry{Con.}

\addplot [color=mycolor2, very thick, mark=x]
  table[row sep=crcr]{%
16	13.9820401\\
36	100.0325611\\
81	976.490821\\
};
\addlegendentry{mld}

\addplot [color=mycolor3, very thick, mark=o]
  table[row sep=crcr]{%
16	11.3111154\\
36	98.6113232000001\\
81	1152.6613124\\
};
\addlegendentry{G}

\addplot [color=mycolor4, very thick, mark=x]
  table[row sep=crcr]{%
16	11.9943983\\
36	79.4314047\\
81	609.8858783\\
196	5600.4618094\\
};
\addlegendentry{R}

\addplot [color=mycolor5, very thick, mark=|]
  table[row sep=crcr]{%
16	5.5773491\\
36	22.1102493\\
81	148.0166021\\
196	988.5957084\\
};
\addlegendentry{M1}

\addplot [color=mycolor6, very thick, mark=o]
  table[row sep=crcr]{%
16	8.9646429\\
36	39.1453961\\
81	224.590377\\
196	1621.9910685\\
};
\addlegendentry{M2}

\addplot [color=mycolor7, very thick, mark=triangle]
  table[row sep=crcr]{%
16	11.0417237\\
36	65.0173671\\
81	471.7219177\\
196	3512.0116667\\
};
\addlegendentry{M3}

\addplot [color=black, dotted]
  table[row sep=crcr]{%
10	4\\
12	6.21608397723432\\
14	8.98564379971739\\
17	14.2239895311328\\
21	23.3239483591063\\
26	38.2608793291309\\
33	66.1464672210808\\
42	114.536469930476\\
54	202.066408145751\\
71	373.287733435714\\
94	697.382308576278\\
127	1357.29000861195\\
175	2747.72160964071\\
200	3681.64799306237\\
};
\addlegendentry{$\text{M}^\text{2}\text{logM}$}

\addplot [color=black, dashed]
  table[row sep=crcr]{%
10	4\\
200	32000\\
};
\addlegendentry{$\text{M}^\text{3}$}

% \addplot [color=black, dotted, forget plot]
%   table[row sep=crcr]{%
% 10	1.33333333333333\\
% 12	2.07202799241144\\
% 14	2.9952145999058\\
% 17	4.74132984371095\\
% 21	7.77464945303545\\
% 26	12.7536264430436\\
% 33	22.0488224070269\\
% 42	38.1788233101586\\
% 54	67.355469381917\\
% 71	124.429244478571\\
% 94	232.460769525426\\
% 127	452.430002870648\\
% 175	915.90720321357\\
% 200	1227.21599768746\\
% };
\end{axis}

\end{tikzpicture}%
    \caption{Construction and solution time for a CSA at \SI{1}{\giga\hertz} over number of array elements with \num{E-3} tolerance for iterative solvers.}
    \label{fig:CSA_scaling2}
\end{figure}

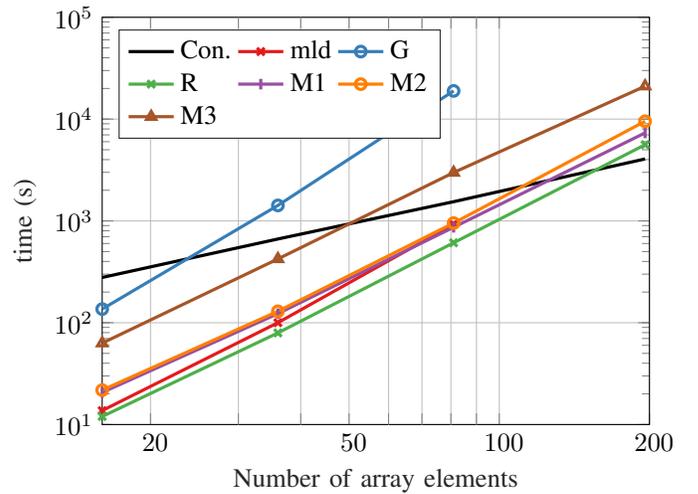
\begin{figure}[t]
    \centering
    % This file was created by matlab2tikz.
%
%The latest updates can be retrieved from
%  http://www.mathworks.com/matlabcentral/fileexchange/22022-matlab2tikz-matlab2tikz
%where you can also make suggestions and rate matlab2tikz.
%
\colorlet{mycolor1}{black}%
\colorlet{mycolor2}{Set1-A}%
\colorlet{mycolor3}{Set1-B}%
\colorlet{mycolor4}{Set1-C}%
\colorlet{mycolor5}{Set1-D}%
\colorlet{mycolor6}{Set1-E}%
\colorlet{mycoloryellow}{Set1-F}%
\colorlet{mycolor7}{Set1-G}%
\begin{tikzpicture}

\begin{axis}[%
width=\columnwidth,
height=7cm,
xmode=log,
xmin=16,
xmax=200,
xminorticks=true,
xlabel style={font=\color{white!15!black}},
xlabel={Number of array elements},
ymode=log,
ymin=10,
ymax=100000,
yminorticks=true,
xmajorgrids,
xminorgrids,
ymajorgrids,
ylabel style={font=\color{white!15!black}},
ylabel={time (s)},
xtick={10,20, 50, 100, 200},
log x ticks with fixed point,
extra x ticks={30, 40, 60, 70, 80, 90},
extra x tick label=\empty,
axis background/.style={fill=white},
legend style={at={(0.03,0.97)}, anchor=north west, legend cell align=left, align=left, draw=white!15!black, legend columns=3}
]
\addplot [color=mycolor1, very thick]
  table[row sep=crcr]{%
16	278.232731\\
36	664.3995005\\
81	1544.2765046\\
196	4060.4189562\\
};
\addlegendentry{Con.}

\addplot [color=mycolor2, very thick, mark=x]
  table[row sep=crcr]{%
16	13.7063739\\
36	100.0493178\\
81	926.3147761\\
};
\addlegendentry{mld}

\addplot [color=mycolor3, very thick, mark=o]
  table[row sep=crcr]{%
16	136.2121309\\
36	1419.4561194\\
81	18909.0858567\\
};
\addlegendentry{G}

\addplot [color=mycolor4, very thick, mark=x]
  table[row sep=crcr]{%
16	11.9943983\\
36	79.4314047\\
81	609.8858783\\
196	5600.4618094\\
};
\addlegendentry{R}

\addplot [color=mycolor5, very thick, mark=|]
  table[row sep=crcr]{%
16	20.6644199\\
36	122.799524\\
81	865.7867649\\
196	7322.2421055\\
};
\addlegendentry{M1}

\addplot [color=mycolor6, very thick, mark=o]
  table[row sep=crcr]{%
16	21.7802024\\
36	129.4978333\\
81	952.3201565\\
196	9519.5469056\\
};
\addlegendentry{M2}

\addplot [color=mycolor7, very thick, mark=triangle]
  table[row sep=crcr]{%
16	62.9961235\\
36	423.1405449\\
81	2984.5173181\\
196	21050.0719914\\
};
\addlegendentry{M3}

% \addplot [color=black, dotted]
%   table[row sep=crcr]{%
% 10	4\\
% 12	6.21608397723432\\
% 14	8.98564379971739\\
% 17	14.2239895311328\\
% 21	23.3239483591063\\
% 26	38.2608793291309\\
% 33	66.1464672210808\\
% 42	114.536469930476\\
% 54	202.066408145751\\
% 71	373.287733435714\\
% 94	697.382308576278\\
% 127	1357.29000861195\\
% 175	2747.72160964071\\
% 200	3681.64799306237\\
% };
% \addlegendentry{$\text{M}^\text{2}\text{logM}$}

% \addplot [color=black]
%   table[row sep=crcr]{%
% 10	4\\
% 200	32000\\
% };
% \addlegendentry{$\text{M}^\text{3}$}

% \addplot [color=black, dotted, forget plot]
%   table[row sep=crcr]{%
% 10	1.33333333333333\\
% 12	2.07202799241144\\
% 14	2.9952145999058\\
% 17	4.74132984371095\\
% 21	7.77464945303545\\
% 26	12.7536264430436\\
% 33	22.0488224070269\\
% 42	38.1788233101586\\
% 54	67.355469381917\\
% 71	124.429244478571\\
% 94	232.460769525426\\
% 127	452.430002870648\\
% 175	915.90720321357\\
% 200	1227.21599768746\\
% };
\end{axis}
\end{tikzpicture}%
    \caption{Construction and solution time for a CSA at \SI{1}{\giga\hertz} over number of array elements with \num{E-6} tolerance for iterative solvers.}
    \label{fig:CSA_scaling}
\end{figure}

\begin{figure}[t]
    \centering
    % This file was created by matlab2tikz.
%
%The latest updates can be retrieved from
%  http://www.mathworks.com/matlabcentral/fileexchange/22022-matlab2tikz-matlab2tikz
%where you can also make suggestions and rate matlab2tikz.
%
\colorlet{mycolor1}{black}%
\colorlet{mycolor2}{Set1-A}%
\colorlet{mycolor3}{Set1-B}%
\colorlet{mycolor4}{Set1-C}%
\colorlet{mycolor5}{Set1-D}%
\colorlet{mycolor6}{Set1-E}%
\colorlet{mycoloryellow}{Set1-F}%
\colorlet{mycolor7}{Set1-G}%
\begin{tikzpicture}

\begin{axis}[%
width=\columnwidth,
height=7cm,
xmode=log,
xmin=16,
xmax=200,
xminorticks=true,
xlabel style={font=\color{white!15!black}},
xlabel={Number of array elements},
ymode=log,
ymin=10,
ymax=100000,
yminorticks=true,
xmajorgrids,
xminorgrids,
ymajorgrids,
ylabel style={font=\color{white!15!black}},
ylabel={time (s)},
xtick={10,20, 50, 100, 200},
log x ticks with fixed point,
extra x ticks={30, 40, 60, 70, 80, 90},
extra x tick label=\empty,
axis background/.style={fill=white},
legend style={at={(0.03,0.97)}, anchor=north west, legend cell align=left, align=left, draw=white!15!black, legend columns=3}
]
\addplot [color=mycolor1, very thick]
  table[row sep=crcr]{%
16	279.8547909\\
25	426.9172895\\
36	621.0734415\\
49	864.252385\\
81	1545.36816\\
121	2295.1132128\\
196	4012.8460936\\
};
\addlegendentry{Con.}

\addplot [color=mycolor2, very thick, mark=x]
  table[row sep=crcr]{%
16	12.9336867\\
25	38.8752841\\
36	103.898247\\
49	242.0872837\\
81	975.3138628\\
};
\addlegendentry{mld}

\addplot [color=mycolor3, very thick, mark=o]
  table[row sep=crcr]{%
16	60.606654\\
25	242.317196\\
36	795.4216458\\
49	2285.8413266\\
};
\addlegendentry{G}

\addplot [color=mycolor4, very thick, mark=x]
  table[row sep=crcr]{%
16	12.685928\\
25	34.5756238\\
36	81.1920725\\
49	173.0261648\\
81	577.8388842\\
121	1656.3120155\\
196	5469.823387\\
};
\addlegendentry{R}

\addplot [color=mycolor5, very thick, mark=|]
  table[row sep=crcr]{%
16	21.1810151\\
25	69.6891721\\
36	208.5957726\\
49	548.591558400001\\
81	2933.9186272\\
121	11993.5319586\\
196	79057.0483424\\
};
\addlegendentry{M1}

\addplot [color=mycolor6, very thick, mark=o]
  table[row sep=crcr]{%
16	22.5957164\\
25	67.2332038\\
36	167.3250551\\
49	381.5178983\\
81	1379.299564\\
121	4777.82784\\
196	27852.2191375\\
};
\addlegendentry{M2}

\addplot [color=mycolor7, very thick, mark=triangle]
  table[row sep=crcr]{%
16	66.3713734\\
25	209.4292986\\
36	620.436253400001\\
49	1355.1256593\\
81	6499.33095670001\\
121	21663.0518864\\
196	291862.6772428\\
};
\addlegendentry{M3}

\addplot [color=black, dotted]
  table[row sep=crcr]{%
10	4\\
12	6.21608397723432\\
14	8.98564379971739\\
17	14.2239895311328\\
21	23.3239483591063\\
26	38.2608793291309\\
33	66.1464672210808\\
42	114.536469930476\\
54	202.066408145751\\
71	373.287733435714\\
94	697.382308576278\\
127	1357.29000861195\\
175	2747.72160964071\\
200	3681.64799306237\\
};
\addlegendentry{$\text{M}^\text{2}\text{logM}$}

\addplot [color=black, dashed]
  table[row sep=crcr]{%
10	4\\
200	32000\\
};
\addlegendentry{$\text{M}^\text{3}$}

% \addplot [color=black, dotted, forget plot]
%   table[row sep=crcr]{%
% 10	1.33333333333333\\
% 12	2.07202799241144\\
% 14	2.9952145999058\\
% 17	4.74132984371095\\
% 21	7.77464945303545\\
% 26	12.7536264430436\\
% 33	22.0488224070269\\
% 42	38.1788233101586\\
% 54	67.355469381917\\
% 71	124.429244478571\\
% 94	232.460769525426\\
% 127	452.430002870648\\
% 175	915.90720321357\\
% 200	1227.21599768746\\
% };
\end{axis}
\end{tikzpicture}%
    \caption{Construction and solution time for a CSA array at \SI{4}{\giga\hertz} over number of array elements with \num{E-3} tolerance for iterative solvers. M1 and M3 did not fully converge for 196 elements.}
    \label{fig:CSA_scaling4}
\end{figure}
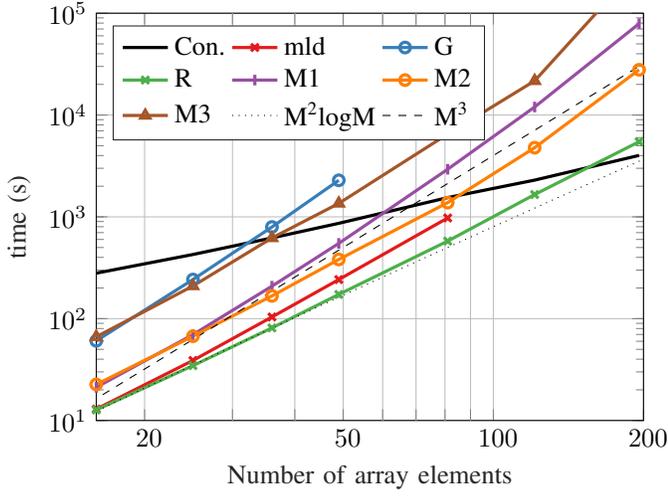

The time to solve the CSA scales differently with regards to the number of array elements depending on the method used. \rOne{At \SI{1}{\giga\hertz}, all accelerated methods have a scaling of around $M^2\log M$ but better than $M^3$ as seen in Fig.~\ref{fig:CSA_scaling2}~and~\ref{fig:CSA_scaling}}, where $M$ is the number of array elements. Though the solution time of the iterative solvers also seems to depend on frequency which the array is simulated at, as seen when comparing the CSA results for \SI{1}{\giga\hertz} in \figref{fig:CSA_scaling2} and \SI{4}{\giga\hertz} in \figref{fig:CSA_scaling4}. For the iterative solvers at \SI{1}{\giga\hertz}, \figref{fig:CSA_scaling2}, the $P_K$ and full RHS variant M1 solves the problem in the least amount of time. The $P_Z$ and full RHS variant M2 takes around $1.6$ times longer time to solve, and the $P_K$ with split RHS variant M3 is the slowest with around $3.6$ times longer time compared to M1. All three MLFFT solvers are faster and scale better compared to regular GMRES. The Rybicki solver has no discernible difference in solution time when comparing \SI{1}{\giga\hertz} and \SI{4}{\giga\hertz}, \rOne{it scales well and is the fastest solver at \SI{4}{\giga\hertz}. As mld and G need to store the full impedance matrix in memory, they can not be evaluated up to $196$ element arrays with the used computer configuration.}

If the tolerance is instead set to \num{E-6}, the iterative results change as seen in \figref{fig:CSA_scaling}. The direct solvers \emph{mldivide} and Rybicki have no defined tolerance and therefore their times are unaffected. For comparisons sake, these and the construction time are still plotted in \figref{fig:CSA_scaling} together with the iterative solvers with \num{E-6} tolerance. The Rybicki method is now the fastest up to $196$ elements and scales similarly as the MLFFT solvers. M1 is still the fastest MLFFT solver, though M2 is now only takes around $1.3$ times longer. Similarly, M3 is also slightly faster and now takes $2.8$ times longer compared to M1. Compared to \figref{fig:CSA_scaling2}, the time for $196$ elements has increased by a factor $7.4$ for M1, $5.9$ for M2, and $6.0$ for M3.

If the frequency is varied rather than the tolerance, there is also a shift in the iterative results as seen when comparing \figref{fig:CSA_scaling2} and \figref{fig:CSA_scaling4}. As the Rybicki method scales slightly worse in these cases, it diverges away from the MLFFT solvers in \figref{fig:CSA_scaling2} and almost intersects M1 in \figref{fig:CSA_scaling}. The direct methods are still largely unaffected. Notably, M2 is faster and clearly scales better compared to M1 as seen in \figref{fig:CSA_scaling4}. Furthermore, M3 scales in a similar manner as in \figref{fig:CSA_scaling2}, but did not converge for each RHS for $196$ elements. Considering that M1 and M3 only differ by wether the RHS is vectorized (M1) or solved one at a time (M3), vectorization seems to have some positive effect on convergence for large arrays.

\begin{figure}[t]
    \centering
    \input{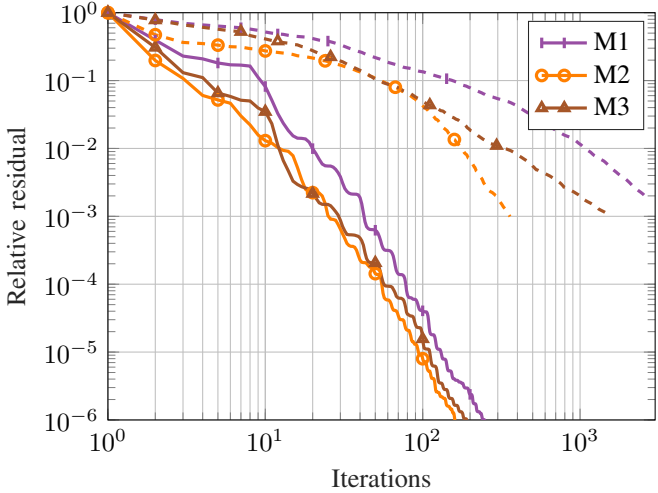}
    \caption{Convergence of MLFFT solvers, relative residual over iteration number for a $196$ element CSA at 1 GHz. For M3, this is the number of interations to solve one RHS. Results for \SI{1}{\giga\hertz} are drawn with solid lines and \SI{4}{\giga\hertz} with dashed lines.}
    \label{fig:CSA_convergence}\label{fig:CSA_convergence4}
\end{figure}

A key performance metric of iterative solvers is their convergence, which is plotted at \SI{1}{\giga\hertz} in \figref{fig:CSA_convergence} for the MLFFT solvers. All solvers converge in a similar manner, though the $P_Z$ preconditioner used in M2 converges faster during the first iterations. Note that while M2 converges in fewer iterations, it takes longer time to apply and is therefore slower overall compared to M1. M3 converges in a similar manner to M1 and is faster to apply to a single RHS, but overall slower due to solving multiple RHS one at a time. At \SI{4}{\giga\hertz}, we see a similar trend in \figref{fig:CSA_convergence4}, although a larger number of iterations is required overall.

\begin{figure}[t]
    \centering
    % This file was created by matlab2tikz.
%
%The latest updates can be retrieved from
%  http://www.mathworks.com/matlabcentral/fileexchange/22022-matlab2tikz-matlab2tikz
%where you can also make suggestions and rate matlab2tikz.
%
\colorlet{mycolor1}{black}%
\colorlet{mycolor2}{Set1-A}%
\colorlet{mycolor3}{Set1-B}%
\colorlet{mycolor4}{Set1-C}%
\colorlet{mycolor5}{Set1-D}%
\colorlet{mycolor6}{Set1-E}%
\colorlet{mycoloryellow}{Set1-F}%
\colorlet{mycolor7}{Set1-G}%
\begin{tikzpicture}

\begin{axis}[%
width=\columnwidth,
height=7cm,
xmode=log,
xmin=16,
xmax=200,
xminorticks=true,
xlabel style={font=\color{white!15!black}},
xlabel={Number of array elements},
ymode=log,
ymin=10,
ymax=200000,
yminorticks=true,
xmajorgrids,
xminorgrids,
ymajorgrids,
xtick={10,20, 50, 100, 200},
log x ticks with fixed point,
extra x ticks={30, 40, 60, 70, 80, 90},
extra x tick label=\empty,
ylabel style={font=\color{white!15!black}},
ylabel={time (s)},
axis background/.style={fill=white},
legend style={at={(0.03,0.97)}, anchor=north west, legend cell align=left, align=left, draw=white!15!black, legend columns=3}
]
\addplot [color=mycolor1, very thick]
  table[row sep=crcr]{%
16	994.043189200001\\
36	2248.5867351\\
81	5264.44502110001\\
196	12974.7806961\\
};
\addlegendentry{Con.}

\addplot [color=mycolor2, very thick, mark=x]
  table[row sep=crcr]{%
16	99.8692134\\
36	923.8142724\\
};
\addlegendentry{mld}

\addplot [color=mycolor3, very thick, mark=o]
  table[row sep=crcr]{%
16	219.3096987\\
36	2672.1783952\\
};
\addlegendentry{G}

\addplot [color=mycolor4, very thick, mark=x]
  table[row sep=crcr]{%
16	103.3367324\\
196	61645.972017\\
};
\addlegendentry{R}

\addplot [color=mycolor5, very thick, mark=|]
  table[row sep=crcr]{%
16	42.1898563\\
36	312.1051621\\
196	27848.527703\\
};
\addlegendentry{M1}

\addplot [color=mycolor6, very thick, mark=o]
  table[row sep=crcr]{%
16	157.5358345\\
36	900.555280800001\\
81	7557.10716580001\\
196	119128.6989791\\
};
\addlegendentry{M2}

\addplot [color=mycolor7, very thick, mark=triangle]
  table[row sep=crcr]{%
16	101.9087959\\
36	837.2207644\\
81	6432.139954\\
196	91794.1874333\\
};
\addlegendentry{M3}

% \addplot [color=black, dotted]
%   table[row sep=crcr]{%
% 10	4\\
% 12	6.21608397723432\\
% 14	8.98564379971739\\
% 17	14.2239895311328\\
% 21	23.3239483591063\\
% 26	38.2608793291309\\
% 33	66.1464672210808\\
% 42	114.536469930476\\
% 54	202.066408145751\\
% 71	373.287733435714\\
% 94	697.382308576278\\
% 127	1357.29000861195\\
% 175	2747.72160964071\\
% 200	3681.64799306237\\
% };
% \addlegendentry{$\text{M}^\text{2}\text{logM}$}

% \addplot [color=black]
%   table[row sep=crcr]{%
% 10	4\\
% 200	32000\\
% };
% \addlegendentry{$\text{M}^\text{3}$}

% \addplot [color=black, dotted, forget plot]
%   table[row sep=crcr]{%
% 10	1.33333333333333\\
% 12	2.07202799241144\\
% 14	2.9952145999058\\
% 17	4.74132984371095\\
% 21	7.77464945303545\\
% 26	12.7536264430436\\
% 33	22.0488224070269\\
% 42	38.1788233101586\\
% 54	67.355469381917\\
% 71	124.429244478571\\
% 94	232.460769525426\\
% 127	452.430002870648\\
% 175	915.90720321357\\
% 200	1227.21599768746\\
% };
\end{axis}

\end{tikzpicture}%
    \caption{Construction and solution time for the wide bowtie array over number of array elements at \SI{150}{\mega\hertz} with \num{E-3} tolerance for iterative solvers. }
    \label{fig:bowtie_scaling}
\end{figure}
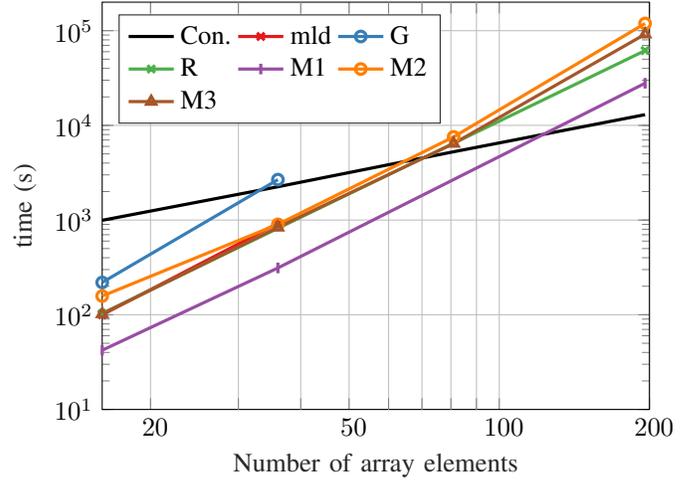

When the tolerance for the CSA is decreased to \num{E-6} in \figref{fig:CSA_scaling}, or when the frequency is increased in \figref{fig:CSA_scaling4}, Rybicki seems to be a better choice. It is faster than M1 by a factor of $0.77$ at $196$ array elements. It does not utilize both levels of the multilevel Toeplitz structure, i.e., storing $Z_A$ scales as $N_xN_y^2N_{e}^2$ for Rybicki compared to $N_xN_yN_{e}^2$ for MLFFT. \rTwo{Here, $N_x$ is the number of columns in the array, 
$N_y$ the number of rows and $N_e$ the number of mesh cells per element.} But when solving the full RHS in MLFFT, the Krylov subspace for $L$ iterations scales as $LN_x^2N_y^2N_{e}$. Therefore, if the number of iterations is on scale with the number of mesh cells of one antenna element, the Krylov subspace starts to dominate the memory usage. On the other hand, the Krylov subspace of the split RHS scales as $LN_xN_yN_{e}$. Therefore, M3 can still be a good choice for memory limited cases, at the cost of computation time in the cases presented in Fig.~\ref{fig:CSA_scaling2},~\ref{fig:CSA_scaling}~and~\ref{fig:bowtie_scaling}.

To verify the results of the Rybicki solver and the MLFFT solvers, the CSA EEP of the center and a corner element are compared with the mld solution. The comparison is done for the $9\times 9$ array, the largest array where the impedance matrix required for the mld solution could be stored in memory. For the EEP from the Rybicki solver, the results are virtually identical, with an RMS deviation on the order of \num{E-14} when compared to the mld solution. For the MLFFT solvers, the RMS deviation is on the order of \num{E-3}, in line with the set convergence criteria. 

Comparing the times for the wide bowtie array ($W=0.375\lambda$) in \figref{fig:bowtie_scaling} with the results from the CSA in \figref{fig:CSA_scaling2} there are some clear differences. The M1 solver is the fastest in both, but for the bowtie array M2, M3 and Rybicki take a similar amount of time to solve. While not presented as a plot, the wide bowtie was also investigated at \SI{225}{\mega\hertz} and similar effects in solution time could be seen as when the CSA frequency was increased to \SI{4}{\giga\hertz}. The three MLFFT solvers see an increase in solution time, \rTwo{making R faster than M2 and M3}. The increased scaling for the MLFFT solvers when going up in frequency was observed for the $P_Z$ preconditioner \rOne{(M2 configuration) whereas the $P_K$ preconditioner (M1 and M3 configuration) had similar scaling for both \SI{150}{\mega\hertz} and \SI{225}{\mega\hertz}. At \SI{225}{\mega\hertz} M1 is still the fastest solver.}

\begin{figure}[t]
    \centering
    % This file was created by matlab2tikz.
%
%The latest updates can be retrieved from
%  http://www.mathworks.com/matlabcentral/fileexchange/22022-matlab2tikz-matlab2tikz
%where you can also make suggestions and rate matlab2tikz.
%
\colorlet{mycolor1}{black}%
\colorlet{mycolor2}{Set1-A}%
\colorlet{mycolor3}{Set1-B}%
\colorlet{mycolor4}{Set1-C}%
\colorlet{mycolor5}{Set1-D}%
\colorlet{mycolor6}{Set1-E}%
\colorlet{mycoloryellow}{Set1-F}%
\colorlet{mycolor7}{Set1-G}%
\begin{tikzpicture}

\begin{axis}[%
width=\columnwidth,
height=5cm,
xmode=log,
xmode=log,
xmin=16,
xmax=1000,
xminorticks=true,
xlabel style={font=\color{white!15!black}},
xlabel={$\text{Number of array elements}$},
xmajorgrids,
xminorgrids,
ymajorgrids,
ymode=log,
ymin=1,
ymax=100000,
yminorticks=true,
ylabel style={font=\color{white!15!black}},
ylabel={time (s)},
axis background/.style={fill=white},
legend style={at={(0.97,0.03)}, anchor=south east, legend cell align=left, align=left, draw=white!15!black, legend columns=3}
]
\addplot [color=mycolor1, very thick]
  table[row sep=crcr]{%
16	100\\
64.0000000000001	393\\
400	2935\\
900.000000000001	9232\\
};
\addlegendentry{Con.}

\addplot [color=mycolor2, very thick, mark=x]
  table[row sep=crcr]{%
16	3.3422502\\
36	21.534735\\
64	104.12991\\
196	2102.8708813\\
};
\addlegendentry{mld}

\addplot [color=mycolor3, very thick, mark=o]
  table[row sep=crcr]{%
16	29.909003\\
36	278.3542488\\
64.0000000000001	1540.4187475\\
196	50008.7018715\\
};
\addlegendentry{G}

\addplot [color=mycolor4, very thick, mark=x]
  table[row sep=crcr]{%
16	3.4667118\\
36	21.8560052\\
64.0000000000001	63.1365337\\
196	952.0991097\\
400	6084.5746261\\
625.000000000001	18117.2783757\\
900.000000000001	44560.2577593\\
};
\addlegendentry{R}

\addplot [color=mycolor5,very thick, mark=|]
  table[row sep=crcr]{%
16	10.8117086\\
36	76.3513432999999\\
64.0000000000001	372.8108777\\
196	7398.94328370001\\
400	50147.5878972\\
};
\addlegendentry{M1}

% \addplot [color=black, dotted]
%   table[row sep=crcr]{%
% 10	4\\
% 12	6.21608397723432\\
% 15	10.5848213315011\\
% 19	18.4652019977589\\
% 25	34.948500216801\\
% 33	66.1464672210808\\
% 44	127.26897526709\\
% 60	256.053780055245\\
% 84	543.108585698106\\
% 121	1219.76242427213\\
% 180	2922.83316661388\\
% 200	3681.64799306237\\
% };
% \addlegendentry{$\text{M}^\text{2}\text{logM}$}

% \addplot [color=black]
%   table[row sep=crcr]{%
% 10	4\\
% 200	32000\\
% };
% \addlegendentry{$M^3$}

% \addplot [color=black, dotted]
%   table[row sep=crcr]{%
% 10	1.33333333333333\\
% 12	2.07202799241144\\
% 15	3.52827377716704\\
% 19	6.15506733258629\\
% 25	11.649500072267\\
% 33	22.0488224070269\\
% 44	42.4229917556968\\
% 60	85.3512600184149\\
% 84	181.036195232702\\
% 121	406.587474757375\\
% 180	974.277722204628\\
% 200	1227.21599768746\\
% };

\end{axis}

\end{tikzpicture}%
    \caption{Construction and solution time for light bowtie array over number of array elements at \SI{150}{\mega\hertz} with \num{E-3} tolerance for iterative solvers. }
    \label{fig:lightBTscaling}
\end{figure}
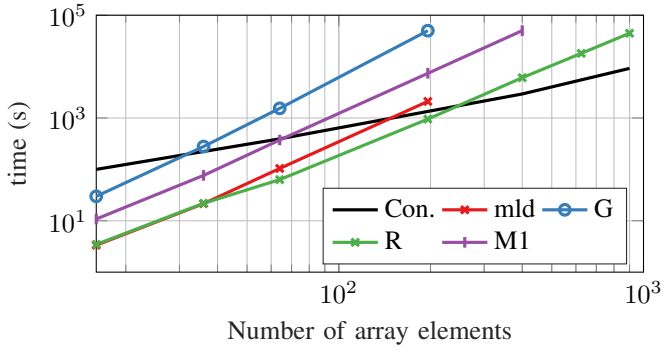

To investigate how the methods work for larger arrays, M1 and R were were used to solve the narrow bowtie up to 400 and 900 elements respectively. As seen in \figref{fig:lightBTscaling} they scale in a similar manner as when solving the CSA or wide bowtie, although Rybicki is notably faster compared to M1.

Which solver to choose for the fastest simulation time seems to depend on the problem. For the cases in \figref{fig:CSA_scaling2} and \figref{fig:bowtie_scaling}, M1 is the fastest investigated solver. Furthermore, the MLFFT solver utilizes both levels in the multilevel Toeplitz structure for the most memory efficient storage of $Z_A$ among the investigated solvers. 

For the iterative solvers, both the smaller $P_K$ and the larger $P_Z$ have their respective advantages and disadvantages. On one hand, $P_K$ is smaller and faster to apply, leading to shorter simulations in \figref{fig:CSA_scaling2} and \figref{fig:bowtie_scaling}. On the other hand, $P_Z$ includes interactions with antenna elements in the same row of the array and allows for around twice as fast convergence in iterations. However, as $P_Z$ is larger it takes longer time to apply to the matrix and is overall slower compared to $P_K$, except for in \figref{fig:CSA_scaling4}. Fewer iterations does decrease the size of the Krylov subspace. A preconditioner which takes more neighboring elements into account, not just those in the same row, may be a better choice compared to both $P_K$ and $P_Z$, though it must also be memory efficient. 

\begin{figure}[t]
    \centering
    \input{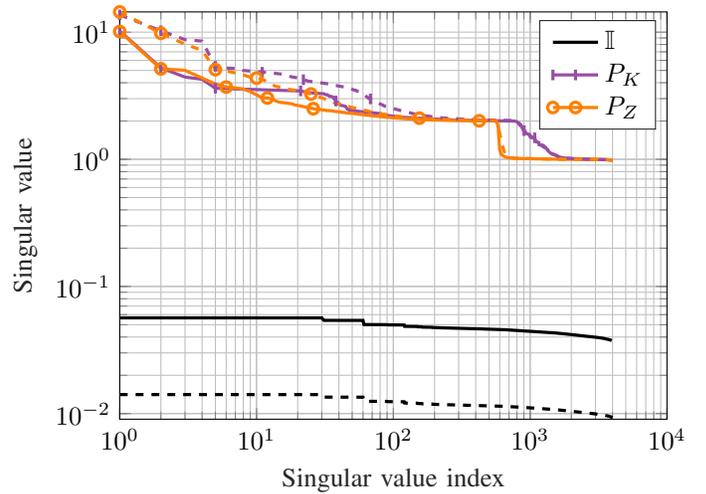}
    \caption{First 4000 singular values of Z for a $6\times 6$ CSA different preconditioners and at different frequencies. Results for \SI{1}{\giga\hertz} are drawn with solid lines and \SI{4}{\giga\hertz} with dashed lines.}
    \label{fig:CSA_SVD}
\end{figure}

To gain further understanding in the difference between the smaller preconditioner $P_K$ and the larger $P_Z$, the spectrum of the preconditioned problem has been investigated. In \figref{fig:CSA_SVD}, the first 4000 of the 17300 singular values of $P^{-1}Z$ are plotted, where $P$ is given by the legend. Here $P=\id$ gives the spectrum of $Z$. These singular values are found using the methods described in \cite{halko2011finding}, implemented for matlab at \cite{Li2021randomized}. The code has been modified to use the fast MLFFT multiplication described in Sec.~\ref{sec:theory}~and~\ref{sec:matinv}. As mentioned in \cite{carson2024towards}, it is difficult to determine convergence behavior from the spectrum alone, but some clear differences and similarities can be seen. The frequency only seems to affect the $\approx 120$ largest singular values. It can also be seen that both preconditioners gives a lot of singular values grouped around 1, which \cite{carson2024towards} gives as a sign of good convergence for GMRES. The preconditioner $P_Z$ has fewer singular values larger than 1, which may explain why it converges in fewer iterations in \figref{fig:CSA_convergence} compared to $P_K$. It can also be seen that the spectrum of $P_Z$ seems to decrease faster than that of $P_K$ from the tenth singular value. Further, it seems that the highest singular values are higher for higher frequency in the preconditioned spectrums. Similar results are seen when investigating the spectrum of the bowtie array.

The three iterative MLFFT configurations are faster for certain problems, mainly when the demands on accuracy are $10^{-3}$ and for lower frequencies. On the other hand, the Rybicki method does seem to be more stable and consistent, with solution times and scaling that is largely unaffected by frequency. Further, it is a direct solver that finds an exact solution down to machine accuracy and that works without any need for a preconditioner.

\section{Conclusion}
\label{sec:conclusion}
We have implemented two variants of solvers that are accelerated by utilizing the multilevel Toeplitz structure inherent to an impedance matrix of an array antenna. The method of \cite{akerstedt2024partitioning} further allows us to apply the solvers to connected elements with an RWG mesh. The goal of these Toeplitz solvers is to find the excitation of arrays with connected elements faster and by using less memory compared to using a full impedance matrix. The results show that by using the decomposition of \cite{akerstedt2024partitioning} and solver utilizing the multilevel Toeplitz structure, we can find the excitation of \rOne{much} larger arrays than conventional regular GMRES or the linear equation solver of MATLAB. \rOne{We also demonstrate that the Topelitz solvers are significantly faster compared to the two conventional solvers.}

We have presented two solvers: the direct Rybicki method and the iterative MLFFT method. The Rybicki method we present is an extension of the \rOne{algorithm} in \cite{press1986numerical} which works on block Toeplitz matrices without requirements on e.g. symmetry. The MLFFT method is presented \rOne{along with the underlying theory, the methods behind it, and ready-to-use code}. The Rybicki method outperforms the direct inversion mldivide in MATLAB \rOne{in memory, time and scaling}. The iterative MLFFT method also outperform regular GMRES in all of its three configurations \rOne{in regards to memory, time and scaling.} For the \rOne{array examples shown}, the best choice between MLFFT configurations and Rybicki depends on array type, frequency, desired accuracy and available memory. The Rybicki method is more robust, delivering high performance regardless of the problem parameters. The MLFFT method has more problem dependent performance, but is the faster choice for certain configurations. Our solver results have been validated against the mldivide results. The problem dependency of the MLFFT solver has been investigated by examining the spectrum of the preconditioned problem.

\appendices
\section{MATLAB Code for Single Level block FFT}
\label{app:SLFFT}
\begin{lstlisting}[style=Matlab-editor]
function Fx = blockFft(x, N_0, fftHandle)
%BLOCKFFT takes the blockwise fft or 
% ifft (1D) with blocksize N_0xN_0
% input:
% x - column block vector
% N_0 - blocksize
% fftHandle - @fft or @ifft
    % allocate Fx
    Fx = x;
    % fft block-wise
    for row = 1:N_0
        Fx(row:N_0:end,:) = fftHandle(x(row:N_0:end,:), [], 1);
    end
end
\end{lstlisting}

\section{MATLAB Code for Dual level block FFT}
\label{app:DLFFT}
\begin{lstlisting}[style=Matlab-editor]
function [Fx] = multiLevelFft(x, N_2, N_1, N_0, fftHandle)
%multiLevelFft takes the multilevel blockwise fft or ifft with
% blocksize N_1N_0xN_1N_0 on the upper level and blocksize N_0xN_0 on the lower
% input:
% x - column block vector
% N_2 - number of blocks on upper level
% N_1 - number of blocks on lower level
% N_0 - size of block on lower level
% fftHandle - @fft or @ifft
    % allocate Fx
    Fx = x;
    % fft blockwise on lower level
    for blockN_2 = 1:N_2
        rowBlockStart = (blockN_2-1)*N_1*N_0;
        for rowInBlock = 1:N_0
            rowIndices = rowBlockStart + (rowInBlock:N_0:((N_1-1)*N_0+rowInBlock));
            Fx(rowIndices,:) = fftHandle(Fx(rowIndices,:), [], 1);
        end
    end
    % fft blockwise on upper level
    N_2Size = N_1*N_0;
    for rowInBlock = 1:N_2Size
        Fx(rowInBlock:N_2Size:end,:) = fftHandle(Fx(rowInBlock:N_2Size:end,:), [], 1);
    end
end
\end{lstlisting}

\bibliographystyle{IEEEtran}
\bibliography{IEEEabrv, sources.bib}

\end{document}